\documentclass[11pt]{article}

\usepackage[cmex10]{amsmath}
\usepackage{amsfonts,amssymb}

\begin{document}

\def\ALERT#1{{\large\bf $\clubsuit$#1$\clubsuit$}}

\numberwithin{equation}{section}
\newtheorem{defin}{Definition}
\newtheorem{theorem}{Theorem}
\newtheorem{proposition}{Proposition}
\newtheorem{notice}{Notice}
\newtheorem{hypothesis}{Hypothesis}
\newtheorem{lemma}{Lemma}
\newtheorem{cor}{Corollary}
\newtheorem{example}{Example}
\newtheorem{remark}{Remark}
\newtheorem{conj}{Conjecture}
\def\begproof{\noindent{\bf Proof: }}
\def\endproof{\par\rightline{\vrule height5pt width5pt depth0pt}\medskip}
\def\div{\nabla\cdot}
\def\rot{\nabla\times}
\def\sign{{\rm sign}}
\def\arsinh{{\rm arsinh}}
\def\arcosh{{\rm arcosh}}
\def\diag{{\rm diag}}
\def\const{{\rm const}}
\def\d{\,\mathrm{d}}
\def\eps{\varepsilon}
\def\phi{\varphi}
\def\theta{\vartheta}
\def\N{\mathbb{N}}
\def\R{\mathbb{R}}
\def\C{\hbox{\rlap{\kern.24em\raise.1ex\hbox
      {\vrule height1.3ex width.9pt}}C}}
\def\P{\hbox{\rlap{I}\kern.16em P}}
\def\Q{\hbox{\rlap{\kern.24em\raise.1ex\hbox
      {\vrule height1.3ex width.9pt}}Q}}
\def\M{\hbox{\rlap{I}\kern.16em\rlap{I}M}}
\def\Z{\hbox{\rlap{Z}\kern.20em Z}}
\def\({\begin{eqnarray}}
\def\){\end{eqnarray}}
\def\[{\begin{eqnarray*}}
\def\]{\end{eqnarray*}}
\def\part#1#2{\frac{\partial #1}{\partial #2}}
\def\partk#1#2#3{\frac{\partial^#3 #1}{\partial #2^#3}} 
\def\mat#1{{D #1\over Dt}}
\def\dx{\nabla_x}
\def\dv{\nabla_v}
\def\grad{\nabla}
\def\Norm#1{\left\| #1 \right\|}
\def\pmb#1{\setbox0=\hbox{$#1$}
  \kern-.025em\copy0\kern-\wd0
  \kern-.05em\copy0\kern-\wd0
  \kern-.025em\raise.0433em\box0 }
\def\bar{\overline}
\def\lbar{\underline}
\def\fref#1{(\ref{#1})}
\def\half{\frac{1}{2}}
\def\oo#1{\frac{1}{#1}}

\def\tot#1#2{\frac{\d #1}{\d #2}} 
\def\laplace{\Delta}
\def\d{\,\mathrm{d}}
\def\N{\mathbb{N}}
\def\R{\mathbb{R}}
\def\supp{\mbox{supp }}
\def\eps{\varepsilon}
\def\phi{\varphi}

\def\J{\mathcal{J}}
\def\I{\mathcal{I}}
\def\K{\mathcal{K}}

\def\Ikn{\I_{k_n}}
\def\sIkn{(\Ikn)_{n\in\N}}

\def\O{\mathcal{O}}
\def\calG{\mathcal{G}}
\def\S{\mathcal{S}}
\def\calW{\mathcal{W}}

\def\td{\alpha}

\def\comment#1{\textbf{#1}}


\centerline{{\huge A Note on the Consensus Finding Problem}}
\centerline{{\huge in Communication Networks}}
\centerline{{\huge with Switching Topologies}}


\vskip 5mm

\centerline{
{\large Jan Haskovec}\footnote{Mathematical and Computer Sciences and Engineering Division,
King Abdullah University of Science and Technology,
Thuwal 23955-6900, Kingdom of Saudi Arabia; 
{\it jan.haskovec@kaust.edu.sa}}
}
\vskip 6mm


\begin{abstract}
In this note, we discuss the problem of consensus finding in communication networks of agents with dynamically switching topologies.
In particular, we consider the case of directed networks with unbalanced matrices of communication rates.
We formulate sufficient conditions for consensus finding in terms of strong connectivity of the underlying directed graphs
and prove that, given these conditions, consensus is found asymptotically.
Moreover, we show that this consensus is an emergent property of the system,
being encoded in its dynamics and not just an invariant of its initial configuration.
\end{abstract}



\section{Introduction}\label{sec:Introduction}
The problem of coordination and consensus finding (``agreement problem'') in distributed communication networks
of dynamic agents has attracted significant attention in many mathematical and engineering communities.
Indeed, apart from being of theoretical interest,
the problem has broad applications in many practical areas,
where groups of agents need to agree upon certain quantities of interest:
cooperative control of unmanned air vehicles, formation control,
flocking and swarming, distributed sensor networks, altitude alignment
of clusters of satellites, congestion control in communication networks and many others.
Consensus problems have a long history in the field of
computer science, particularly in automata theory and distributed computation. 
The critical problem for coordinated control is to design appropriate
protocols and algorithms such that the group of agents can reach consensus
on the shared information in the presence of limited and/or unreliable information exchange
and dynamically changing interaction topologies.
Therefore, it is important to address agreement problems in
their general form, with networks of dynamic agents with directed
information flow under possible link failure and creation (i.e., switching
network topology). Consensus problems in various contexts
have recently been addressed in \cite{B1}--\cite{Vicsek}, to name a few.

The topology of the interaction (information exchange via communication or direct sensing)
between the agents is typically represented by directed graphs.
Clearly, finding a global network consensus by any means
is only possible if the communication graph is in some sense ``sufficiently connected''.
Loosely speaking, each node must possess, over some
sufficiently dense collection of time intervals,
a communication path to every other node in the network. 
A preliminary result for consensus finding was presented in \cite{B8}, where
a linear update scheme was proposed for directed graphs. However, the
analysis in \cite{B8} was not able to utilize all available communication links.
A solution to this issue was presented in \cite{B4} for time-invariant communication
topologies. Information consensus for dynamically evolving
information was addressed in \cite{B9} in the context of spacecraft formation
flying where the exchanged information is the configuration of the virtual
structure associated with the (dynamically evolving) formation.
In many applications, the interaction topology between agents may
change dynamically. For example, communication links
may be unreliable due to disturbances and/or subject to communication
range limitations. If information is being exchanged by
direct sensing, the locally visible neighbors will likely
change over time, typically due to limited sensing radius or occlusion.
This is usually the case in models of biological flocking and swarming,
see, e.g. the recent surveys \cite{Carrillo-review, Vicsek-survey}.
A well known example is the Vicsek model \cite{Vicsek}
with possible changes over time in each agent's nearest neighbors.
A theoretical explanation for its behavior was provided in \cite{B2},
where it was shown that consensus can be achieved if the union of the interaction graphs
for the team is connected frequently enough as the system evolves.
Another well known model of biological collective interaction is the
Cucker-Smale model, \cite{CS1, CS2}, originally formulated as a model
for language evolution, but later established as a flocking model.
Its recent modification \cite{Albi-Pareschi, Haskovec} introduces topological interactions
between the agents and transforms it into a communication network with switching topologies.
In \cite{B15} and \cite{Ren-Beard05}, the problem of information consensus finding
among multiple agents with dynamically changing directed interaction topologies
was considered. It was shown there that consensus is achieved asymptotically if the union
of the directed interaction graphs have a spanning tree frequently enough
as the system evolves. The proof was based on algebraic graph theory
and abstract theory of stochastic matrices.

In this note, we provide a relatively simple proof of asymptotic consensus finding
in directed communication networks with dynamically switching topology. Our main structural assumption
is that the directed graphs representing the communication topology
are strongly connected frequently enough. This is a stronger assumption
than the existence of a spanning tree, as made in \cite{B15, Ren-Beard05}. However, our proof
is much simpler than that of \cite{B15, Ren-Beard05}, being based merely on tools of elementary calculus.
Moreover, in \cite{B15, Ren-Beard05} it is required that the time intervals where
a particular interaction graph governs the communication are bounded; we do not pose
such a requirement.

\section{Problem statement}
We consider a system of $N\in\N$ agents, each agent being in state $\xi_i\in \R$, $i=1,\dots,N$.
The agents dynamically exchange information about their states, and the interaction/communication
topology at time $t\ge 0$ is described by a directed graph $G(t)$, consisting of $N$ nodes.
In $G(t)$, the $i$-th node $A_i$ represents the $i$-th agent and a directed edge from $A_i$
to $A_j$ represents a unidirectional information exchange from $A_i$ to $A_j$,
i.e., agent $j$ is as time $t$ receiving information from agent $i$.
Moreover, each edge carries a (time dependent) weight $g_{ij}(t)$,
representing the relative intensity of the information exchange between the agents $i$ and $j$;
we set $g_{ij}(t)=0$ if there is no communication between $i$ and $j$ (and, thus, no edge between
$A_i$ and $A_j$ in the graph $G(t)$) at time $t\geq 0$.
Clearly, with $N$ agents we can have at most $N!$ different topologies,
which we denote by $G^k$, $k=1,\dots,N!$,
and the network switches between (some of) them during its temporal evolution.

A directed path in graph $G$ is a sequence of edges $(A_{i_1},A_{i_2})$,
$(A_{i_2},A_{i_3})$, $\dots$ in that graph. The graph is called
strongly connected if for any pair of distinct nodes $A_i$, $A_j$ there is a directed path 
from $A_i$ to $A_j$ and a directed path from $A_j$ to $A_i$.
A directed tree is a directed graph, where every node,
except the root, has exactly one parent. A spanning tree of a directed
graph is a directed tree formed by graph edges that connect all the
nodes of the graph. We say that a graph has (or contains)
a spanning tree if a subset of the edges forms a spanning tree.

In our paper, we will consider the following continuous-time
communication scheme, which is widely used in the context
of consensus seeking,
\[
   \dot \xi_i(t) = \frac{1}{\sum_{j=1}^N g_{ij}(t)} \sum_{j=1}^N g_{ij}(t) (\xi_j(t)-\xi_i(t)),
\]
where all the functions $g_{ij}: [0,\infty) \to [0,\infty)$ are piecewise
continuous. In fact, by a rescaling of the functions $g_{ij}$,
we can rewrite the above system as
\( \label{scheme}
   \dot \xi_i(t) = \sum_{j=1}^N g_{ij}(t) (\xi_j(t)-\xi_i(t)),
\)
with the normalization
\( \label{normalization}
   \sum_{j=1}^N g_{ij}(t) = 1 \qquad\mbox{for all } t\geq 0\mbox{ and } i=1,\dots,N,
\)
which is the form we will use in the sequel.
The agent system is said to achieve \emph{asymptotic consensus}
if for any initial datum $(\xi_1(0),\dots,\xi_N(0))\in\R^N$,
there is a $\xi^\infty\in\R$ such that
\(  \label{consensus}
    \lim_{t\to\infty} \xi_i(t) = \xi^\infty,\qquad i=1,\dots,N.
\)

Since the system \eqref{scheme} is in general a system of ordinary
differential equations with discontinuous coefficients, we shall
state a precise definition of its solution:

\begin{defin}\label{def:sol}
Denote by $\{\J_k\}_{k\in\N}$ the at most countable system of open intervals
such that all $g_{ij}$ in \eqref{scheme} are continuous on every $\J_k$ and $\bigcup_{k\in\N} \overline{\J_k} = [0,\infty)$.
We call the globally continuous curve $(\xi_1(t), \dots, \xi_N(t))_{t\geq 0} \subset \R^{N}$ a solution to \eqref{scheme}
if it solves the ODE system on every open interval $\J_k$, $k\in\N$.
\end{defin}

Moreover, for all $t\geq 0$, let us denote by $G(t)$ the directed graph induced by the edge weights $g_{ij}(t)$,
i.e., there is a directed edge connecting vertex $A_i$ to vertex $A_j$ at time $t$ if and only if $g_{ij}(t) > 0$.

\section{Main result}
\begin{theorem}\label{thm:flocking}
Let all $g_{ij}(t)$ be piecewise continuous nonnegative functions.
Let $(\xi_1(t), \dots, \xi_N(t))_{t\geq 0}\subset \R^{N}$ denote the solution of the system \eqref{scheme}
in the sense of Definition \ref{def:sol},
subject to the initial condition $(\xi_1(0),\dots,\xi_N(0))\in\R^N$.
Assume that there exists a topological configuration with a strongly connected directed graph,
say $G^0$, where the system spends an infinite amount of time, i.e.,
\[
    \mbox{meas}\{t\geq 0;\; G(t) = G^0 \} = +\infty,
\]
and, moreover, $G^0$ has the property that each nonzero weight $g_{ij}(t)$ satisfies $g_{ij}(t) \geq c_{ij} > 0$
for all $t$ such that $G(t) = G^0$.

Then the system reaches an asymptotic velocity consensus, i.e., there exists a $\xi^\infty\in\R$
in the convex hull of $\{\xi_1(0),\dots,\xi_N(0)\}$ such that
\(   \label{thm1:statement}
    \lim_{t\to\infty} \xi_i(t) = \xi^\infty\qquad \mbox{for all }  i=1,\dots,N\,.
\)
\end{theorem}

\begproof
The proof will be carried out in three steps.

\textbf{Step 1: Maximal value.}
Due to the assumed piecewise continuity of the functions $g_{ij}$,
we can construct an at most countable system of disjoint open time intervals
$\I_k := (t_{k-1}, t_k)$, $k\in\N$
where the system does not change its topological configuration,
i.e., $G(t) \equiv G_k$ on $\I_k$ for some fixed directed graph $G_k$.
Moreover, we may assume that $\bigcup_{k\in\N} \overline{\I_k} = [0,\infty)$.
Inspired by \cite{CFRT}, let us define for all $t\geq 0$ the function
\[
    \omega(t) := \max_{i=1,\dots,N} |\xi_i(t)|\,,
\]
and, moreover, denote
\[
  M(t) := \mbox{argmax}_{i=1,\dots,N} |\xi_i(t)|\,.
\]
If $M(t)$ is not uniquely determined
(i.e., there are several maximal values $|\xi_i(t)|$), we choose one of the indices arbitrarily,
but in such a way that $M(t)$ stays constant on the longest time interval.
Since the number of agents $N$ is finite and the curves $\xi_i(t)$ are continuous,
there exists an at most countable system of open disjoint intervals $(\K_k)_{k\in\N}$
such that $\bigcup_{k\in\N} \overline{\K_k} = [0,\infty)$,
and $M(t)$ is constant on every $\K_k$. 
To ease the notation, we will usually skip the explicit dependence of $M$ on $t$ (or $k$) in the sequel.

By intertwining the two systems $(\I_k)_{k\in\N}$ and $(\K_k)_{k\in\N}$, we construct another at most countable system of disjoint intervals,
denoted by abuse of notation again by $(\I_k)_{k\in\N}$, such that the topological configuration does not change
\emph{and} the index $M$ is constant on each $\I_k$.
Then, on every $\I_k$ we have
\[
    \frac12 \tot{}{t} \omega(t)^2 = \frac12 \tot{}{t} |\xi_M|^2 &=& \sum_{j=1}^N g_{Mj} (\xi_j - \xi_M)\xi_M \\
	  &\leq& \sum_{j=1}^N g_{Mj} \left( |\xi_j| - |\xi_M| \right) |\xi_M| \,. 
\]
Dividing by $|\xi_M|$ (note that if $|\xi_M|$ was zero, there would be nothing to prove), we obtain
\[
    \tot{}{t} \omega(t) = \tot{}{t} |\xi_M| \leq \sum_{j=1}^N g_{Mj} \left( |\xi_j| - |\xi_M| \right) \leq 0
\]
on every $\I_k$, where the nonpositivity of the right-hand side is due to the maximality of $|\xi_M|$.
Observe that the above inequality holds universally, regardless of whether the configuration
is (strongly) connected or not.
Consequently, $\omega(t)$ is a globally continuous, nonincreasing and nonnegative function,
so that there exists an $0 \leq \omega_\infty \leq \omega(0)$ such that
$\lim_{t\to\infty} \omega(t) = \omega_\infty$.

\textbf{Step 2: Induction.}
Since, by assumption, the system spends an infinite amount of time in the strongly connected configuration $G^0$,
we can pick the corresponding subsystem out of $(\I_k)_{k\in\N}$, with infinite length. By a further
subselection we get the system $\sIkn$ of infinite length, where $M(t) \equiv M_0$ for some fixed $1\leq M_0 \leq N$.
Therefore, denoting $\I^0 := \bigcup_{n\in\N}\Ikn$, we have the configuration $G^0$ and the maximal value index $M_0$ for all $t\in\I^0$.
Moreover,
\[
    |\xi_{M_0}(t)| \to \omega_\infty \qquad\mbox{and}\qquad \tot{}{t} |\xi_{M_0}(t)|^2 \to 0
\]
as $t\to\infty,\; t\in\I^0$,
where the convergence of the time derivative is due to the monotonicity of $|\xi_{M_0}(t)|$.
Since the configuration $G^0$ is fixed and strongly connected,
there exists an index $j_0$ such that $g_{M_0,j_0} =: g^0 >0$ on $\I^0$.
Then we have, for $t\in\I^0$,
\[
    \frac12 \tot{}{t} |\xi_{M_0}(t)|^2 &=& \sum_{j=1}^N g_{M_0,j} (\xi_j-\xi_{M_0}) \xi_{M_0}(t)   \\
    &\leq& g^0 (\xi_{j_0}-\xi_{M_0}) \xi_{M_0}(t) \leq 0,
\]
where we used the inequality $(\xi_j-\xi_{M_0})\cdot \xi_{M_0} \leq 0$ implied by the maximality of $|\xi_{M_0}|$.
Now, since the left-hand side tends to zero as $t\to\infty$ and $g^0$ is bounded from below by a positive constant, we have
\[
    \xi_{j_0} \xi_{M_0}(t) - \xi_{M_0}(t)^2 \to 0 \qquad\mbox{as } t\to\infty,\; t\in\I^0,
\]
and this further implies $\xi_{j_0} \xi_{M_0}(t) \to \omega_\infty^2$.
Finally, since $\xi_{j_0}\xi_{M_0} = \omega_\infty^2$ if and only if $\xi_{j_0} = \xi_{M_0}$,
we have
\(   \label{convergence one velocity}
    (\xi_{j_0} - \xi_{M_0})(t) \to 0 \qquad\mbox{as } t\to\infty,\; t\in\I^0.
\)
Moreover, we calculate
\[
    \frac12 \tot{}{t} \xi_{j_0}(t)^2 &=& \sum_{l=1}^N g_{j_0,l} (\xi_l-\xi_{j_0}) \xi_{j_0}(t) \\
           &=& \sum_{l=1}^N g_{j_0,l} \bigl[ (\xi_l-\xi_{M_0})\xi_{M_0} \\&& +\; (\xi_l-\xi_{M_0})(\xi_{j_0}-\xi_{M_0}) + (\xi_{M_0}-\xi_{j_0})\xi_{j_0} \bigr] \\
           &\leq& \sum_{l=1}^N g_{j_0,l} (\xi_l-\xi_{M_0})\xi_{M_0}(t) + 3\omega(0)|\xi_{j_0}-\xi_{M_0}| \\
           &\leq& 3\omega(0)|\xi_{j_0}-\xi_{M_0}|,
\]
where we used the estimate $\max\left(|\xi_{j_0}|,|\xi_l|\right) \leq \omega(0)$ and the maximality of $\xi_{M_0}$.
By \eqref{convergence one velocity} we have then
\[
   \limsup_{t\to\infty,\,t\in\I^0} \tot{}{t} \xi_{j_0}(t)^2 \leq 0. 
\]
Since $|\xi_{j_0}(t)| \to \omega_\infty$ from below on $\I^0$ as $t\to\infty$, we conclude that
\(   \label{convergence derivative}
     \lim_{t\to\infty,\,t\in\I^1} \tot{}{t} \xi_{j_0}(t)^2 = 0.
\)
where $\I^1$ is a system of subintervals of $\I^0$, still of infinite Lebesgue measure.

We will now show that \eqref{convergence one velocity} in fact holds for any index $\hat\jmath\in\{1,\dots,N\}$.
Due to the simple connectivity of the digraph $G^0$, there exists a path
$M_0 \mapsto j_0 \mapsto j_1 \mapsto \dots \mapsto j_\ell \mapsto \hat\jmath$, such that
$g_{M_0,j_0}>0$, $g_{j_0,j_1}>0$, $\dots$, $g_{j_\ell, \hat\jmath}>0$ on $\I^0$.
We proceed inductively, showing first that the results \eqref{convergence one velocity} and \eqref{convergence derivative} hold for $j_1$ as well.
Indeed, passing to the limit in 
\[
    \frac12 \tot{}{t} \xi_{j_0}(t)^2 &\leq& \sum_{l=1}^N g_{j_0,l} (\xi_l-\xi_{M_0})\xi_{M_0}(t) \\&& +\; 3\omega(0)|\xi_{j_0}-\xi_{M_0}| \\ 
       &\leq& g_{j_0,j_1} (\xi_{j_1}-\xi_{M_0})\xi_{M_0}(t) \\&&\; + 3\omega(0)|\xi_{j_0}-\xi_{M_0}|,
\]
we obtain, due to \eqref{convergence derivative} and the maximality of $\xi_{M_0}$,
\[
    0 \leq \lim_{t\to\infty,\,t\in\I^1} (\xi_{j_1}-\xi_{M_0})\xi_{M_0}(t) \leq 0,
\]
which immediately gives \eqref{convergence one velocity} for $\xi_{j_1}$ on $\I^1$.
Using this result, we argue as before to conclude
\[
   \lim_{t\to\infty,\,t\in\I^2} \tot{}{t} \xi_{j_1}(t)^2 = 0
\]
with $\I^2$ a system of subintervals of $\I^1$ of infinite Lebesgue measure.
This is \eqref{convergence derivative} for $\xi_{j_1}$ on $\I^2$.
Proceeding inductively, after a finite number of steps we reach the index $\hat\jmath$.

We conclude that there exists a sequence $(t_k)_{k\in\N}\subset \I^0$, $t_k\to\infty$, such that, for all $j=1,\dots,N$,
\(  \label{almost_there}
    \xi_j(t_k) - \xi_{M_0}(t_k) \to 0 \qquad\mbox{as } k\to\infty. 
\)

\textbf{Step 3: Conclusion.}
The fact that $\sum_{j=1}^N g_{ij} xi_j$ is a convex combination of the values $\xi_1,\dots,\xi_N$
directly implies that the convex hull of $\{\xi_1,\dots,\xi_N\}$ is nonexpanding in time,
\[
   \mbox{ch}\{\xi_1,\dots,\xi_N\}(t) \subseteq \mbox{ch}\{\xi_1,\dots,\xi_N\}(s) \qquad\mbox{for all } t>s\geq 0.
\]
Due to \eqref{almost_there}, its diameter shrinks to zero as $t\to\infty$,
and, consequently, there exists a value $\xi^\infty\in\mbox{ch}\{\xi_1,\dots,\xi_N\}(0)$ such that
\[
    \lim_{t\to\infty} \xi_j(t) = \xi^\infty \qquad\mbox{for all } j=1,\dots,N. 
\]
\endproof

Let us note that the asymptotic consensus value $\xi^\infty$ in the above Theorem cannot be,
in general, predicted from the initial datum $(\xi_1(0),\dots,\xi_N(0))$,
beyond the trivial fact that $\xi^\infty$ is a convex combination of $\xi_i(0)$, $i=1,\dots,N$.
We may therefore consider $\xi^\infty$ as an emergent property
of the communication network, in the sense that the asymptotic consensus is encoded in the
dynamics of the network and not just as an invariant of its initial configuration.
The existence of invariants of the communication scheme \eqref{scheme} is related
to the balance of the matrix $(g_{ij}(t))_{i,j=1}^N$ of the communication rates.
Indeed, a quick calculation reveals that
\[
   \sum_{i=1}^N \dot\xi_i(t) = \sum_{j=1}^N \left(\sum_{i=1}^N g_{ij} - \sum_{i=1}^N g_{ji}\right)\xi_j
\]
and this is always zero if and only if the matrix $(g_{ij}(t))_{i,j=1}^N$ is balanced, i.e.,
its row and column sums are the same. If this is the case for all $t\geq 0$, then the average 
$\frac{1}{N} \sum_{i=1}^N \xi_i(t)$ is an invariant of the evolution, and we have
\[
    \xi^\infty = \frac{1}{N} \sum_{i=1}^N \xi_i(0)\,,
\]
as for instance in the classical Cucker-Smale model \cite{CS1, CS2},
where the communication matrix is symmetric, thus balanced.


\end{document}